\documentclass{amsart}
\usepackage{amsmath}
\usepackage{amssymb}
\usepackage{eucal}
\usepackage[all]{xy}

\usepackage[OT2,OT1]{fontenc}
\DeclareMathAlphabet{\cyr}{OT2}{wncyr}{m}{n}

\newcommand{\ZZ}{\mathbb{Z}}
\newcommand{\FF}{\mathbb{F}}
\newcommand{\CC}{\mathbb{C}}
\newcommand{\QQ}{\mathbb{Q}}

\newcommand{\GG}{\mathbb{G}}
\newcommand{\TT}{\mathbb{T}}

\newcommand{\cA}{\mathcal{A}}
\newcommand{\cB}{\mathcal{B}}
\newcommand{\cO}{\mathcal{O}}

\newcommand{\fp}{\mathfrak{p}}
\newcommand{\fq}{\mathfrak{q}}

\newcommand{\Sha}{\cyr{X}}

\newcommand{\tA}{\widetilde{A}}
\newcommand{\tB}{\widetilde{B}}
\newcommand{\tC}{\widetilde{C}}
\newcommand{\tJ}{\widetilde{J}}

\newcommand{\af}{\hat{f}}
\newcommand{\ag}{\hat{g}}

\DeclareMathOperator{\coker}{coker}
\DeclareMathOperator{\End}{End}
\DeclareMathOperator{\Gal}{Gal}

\DeclareMathOperator{\Hom}{Hom}
\DeclareMathOperator{\Ext}{Ext}

\DeclareMathOperator{\Spec}{Spec}
\DeclareMathOperator{\Sym}{Sym}

\newcommand{\iso}{\stackrel{\sim}{\to}}

\newcommand{\Qbar}{\overline{\QQ}}
\newcommand{\Kbar}{\overline{K}}

\newcommand{\power}[2]{#1[\![ #2 ]\!]}

\theoremstyle{plain}
\newtheorem{theorem}{Theorem}
\newtheorem{lemma}[theorem]{Lemma}
\newtheorem{corollary}[theorem]{Corollary}
\newtheorem{proposition}[theorem]{Proposition}

\theoremstyle{definition}

\theoremstyle{remark}
\newtheorem*{remark}{Remark}
\newtheorem*{ack}{Acknowledgements} 
\newtheorem*{example}{Example}

\begin{document}


\title[Extensions of abelian varieties]{Extensions of abelian varieties
  defined over a number field}
\author{Matthew A.~Papanikolas}
\address{Department of Mathematics\\
  Texas A{\&}M University\\
  College Station, TX 77843}
\email{map@math.tamu.edu}
\author{Niranjan Ramachandran}
\address{Department of Mathematics\\
  University of Maryland\\
  College Park, MD 20742}
\address{Max-Planck-Instit\"ut f\"ur Mathematik\\
  Vivatsgasse 7\\
  D-53111 Bonn, Germany} 
\email{atma@math.umd.edu}
\thanks{Research of the first author partially supported by NSF and NSA
  grants.  Research of the second author partially supported by MPIM}

\date{April 20, 2004; August 26, 2004 (final)}

\begin{abstract}
  We study the arithmetic aspects of the finite group of extensions of
  abelian varieties defined over a number field. In particular, we
  establish relations with congruences between modular forms and
  special values of $L$-functions.
\end{abstract}


\maketitle

\section{Introduction} 

J.~Milne \cite{milne68, milne68ts} has established striking
connections between arithmetic and the extensions of abelian varieties
over finite fields (see Theorem~\ref{Milne}).  Our aim here is to
relate extensions of abelian varieties over number fields to
congruences between modular forms and special values of
\(L\)-functions of motives.

For abelian varieties $A$ and $B$ over a field $K$, consider the group
$\Ext^1_K(A,B)$ of isomorphism classes of Yoneda extensions of $A$ by
$B$ in the category of commutative group schemes over~$K$.  This group
is torsion by the Poincar\'e reducibility theorem, but it need not be
finite.  For instance, $\Ext^1_{\CC}(A,B) \cong ({\QQ}/{\ZZ})^4$ when
$A$ and $B$ are non-isogenous elliptic curves over $\CC$.  If $A$ is
an elliptic curve with complex multiplication over $\CC$, then a
beautiful result of S.~Lichtenbaum \cite[Thm.~6.1]{pr:WB} states that
$\Ext^1_{\CC}(A,A)$ is naturally isomorphic to the torsion subgroup of
$A(\CC)$.

When \(K\) is a number field, the group $\Ext^1_K(A,B)$ is finite
(Theorem \ref{T:ExtNF}).  We show that the order of 
$\Ext^1_K(A,B)$ is related to the following objects:
\begin{itemize}
\item  (Corollary \ref{C:Extmodf}) congruences between modular forms,
  when \(A\) and \(B\) are elliptic  curves over \(\QQ\);

\item (Theorems \ref{T:ParExtRestMod}, \ref{T:ParExtTMod}) the congruence
  modulus and modular degree of an elliptic curve  over \(\QQ\); 

\item (Theorem \ref{T:Sym2}) the special value of the \(L\)-function
  \(L(\Sym^2 E,s)\) at \(s=2\) for certain elliptic curves \(E\) over
  \(\QQ\).  
\end{itemize}
It is important to note that the third result is simply a restatement
of the deep results of F.~Diamond, M.~Flach, and L.~Guo \cite{dfg01,
dfg} in our context.

\begin{ack}
  We thank J.~Milne for graciously allowing the unpublished proof of
  Theorem \ref{T:ExtNF} in \S\ref{S:pre} to be included here.  We
  thank K.~Ribet for sharing his proof of Corollary~\ref{C:Extmodf} in
  \S\ref{S:ExtCong} with us and for permitting us to include it to
  supplement our own.  We also thank S.~Ahlgren, V.~Balanzew,
  N.~Dummigan, S.~Lichtenbaum, G.~Pappas, and J.~Teitelbaum for their
  helpful remarks and suggestions. We thank the referee for many
  useful suggestions.
\end{ack}

\section{Preliminaries} \label{S:pre}

\subsection*{Extensions over finite fields} 

We recall the fundamental results of Milne \cite{milne68, milne68ts}
on \(\Ext(A,B)\) over a finite field.

\begin{theorem}[Milne (a) {\cite[Thm.~3]{milne68}}; (b)
  {\cite{milne68ts}}]\label{Milne}  Let $A$ and $B$ be abelian
  varieties over a finite field $\FF_q$.
\begin{enumerate}
\item[(a)] The group \(\Ext^1_{\FF_q}(A,B)\)  is
  finite, and its order is given by 
\[
\#\Ext^1_{\FF_q}(A,B)= \pm \frac{q^{d_Ad_B}}{D} \prod_{a_i
  \neq b_j} \left(1- \frac{a_i}{b_j}\right),
\]
where \(a_i\) and \(b_j\) are the eigenvalues of Frobenius associated
with \(A\) and \(B\) over \(\FF_q\), \(d_A\) and \(d_B\) are the
dimensions of $A$ and $B$, and \(D\) is the discriminant of the trace
pairing \(\Hom_{\FF_q}(A,B) \times \Hom_{\FF_q}(B,A) \to \End(A) \to
\ZZ\).
\item[(b)] Let \(J\) be the Jacobian of a smooth projective curve
  \(C\) over \(\FF_q\).  The group \(\Ext^1_{\FF_q}(J,A)\) is
  isomorphic to the Tate-Shafarevich group \(\Sha(A/{\FF_q(C)})\) of
  the constant abelian variety \(A\) over the function field
  \(\FF_q(C)\) of \(C\).
\end{enumerate}
\end{theorem}

\begin{remark}
  Every abelian variety \(A\) defines integral motives \(h^1(A)\) and
  \(h_1(A)\); note that \(h_1(A)\) is isomorphic to
  \(h^1(A^{\vee})(1)\) defined by the dual abelian variety
  \(A^{\vee}\). Part (a) of Theorem \ref{Milne} relates the order of 
\(\Ext^1_{\FF_q}(A,B)\) to the special value at \(s=0\) of the
  \(L\)-function \(L(h_1(A)\otimes h^1(B),s)\) of the tensor product
  motive \cite[Thm.~10.1]{milram}.
\end{remark}
 
\subsection*{Finiteness of $\Ext^1_K(A,B)$ over number fields}

For any abelian group (scheme) $G$, let $G_n$ denote the kernel of
multiplication by the integer $n$.  For a prime $p$, let $T_pG =
\varprojlim G_{p^m}$ denote the associated Tate module, and let $TG =
\varprojlim G_n$ denote the total Tate module.

 \begin{theorem}[Milne-Ramachandran] \label{T:ExtNF}
  Let $A$ and $B$ be abelian varieties over a number field $K$.  Then
  $\Ext^1_K(A,B)$ is finite.
\end{theorem}
\begin{proof}
For any integer $n$, taking $\Ext^i_K(-,B)$ of the Kummer sequence
\[
  0 \to A_n \to A \stackrel{n}{\to} A \to 0
\]
gives a short exact sequence 
\begin{equation} \label{E:ExtABn}
  0 \to \Hom_K(A,B) \otimes \ZZ/n\ZZ \stackrel{\alpha_n}{\to} \Hom_K(A_n,B_n) \to
  \Ext^1_K(A,B)_n \to 0,
\end{equation}
using $\Hom_K(A_n,B) = \Hom_K(A_n,B_n)$.
Now taking the inverse limit of \eqref{E:ExtABn} over powers $n= p^m$
of a prime $p$ and using that $\Hom_{K}(A_{p^m},B_{p^m}) =
\Hom_{K}(T_p A, B_{p^m})$, we have an exact sequence
\begin{equation} \label{E:TpExt}
  0 \to \Hom_K(A,B) \otimes \ZZ_p \stackrel{\alpha}{\to} \Hom_K(T_p A,T_p B) \to T_p
  \Ext^1_K(A,B) \to 0.
\end{equation}By 
Faltings' theorem \cite[\S IV.1]{falt}, \(\alpha\) is an isomorphism,
and so $T_p \Ext^1_K(A,B) = 0$ for all $p$.  As $\Ext^1_K(A,B)_n$ 
is finite for all $n$ by \eqref{E:ExtABn}, the $p$-primary subgroup
\(\Ext^1_K(A,B)(p)\) of $\Ext^1_K(A,B)$ is finite for all $p$.  To
prove the theorem, it now suffices to show that \(\Ext^1_K(A,B)(p)\)
is nonzero only for finitely many \(p\).  This follows from \cite[\S
IV.4]{falt} which says that \(\alpha_p\) in \eqref{E:ExtABn} is an
isomorphism for sufficiently large primes \(p\).
\end{proof}

\begin{remark} 
  (i) In the case of an elliptic curve $A$ over $\QQ$, the results of
  J.-P.~Serre \cite{serre72} provide upper bounds on the order of
  $\Ext^1_{\QQ}(A,A)$.  (ii) There are explicit bounds (see
  \cite[Cor.~1]{maswus95}) for the primes dividing $\#\Ext^1_K(A,B)$.
  (iii) The group \(\Ext^1_K(A,B)\) is not invariant under isogeny.
\end{remark}

\subsection*{Representability issues}

Fix abelian varieties \(A\) and \(B\) over a field \(K\).  The group
\(\Ext^1_K(A,B)\) is similar to the group \(\Ext^1_K(A,\GG_m)\), but
there are some differences.  For example, the functor \(S \mapsto
\Ext^1_S(A\times S, \GG_m)\) on the category of schemes over \(K\) is
representable (by the dual abelian variety).  However, the
corresponding functor for \(\Ext^1_K(A,B)\) is not representable.  If
\(L\) is an extension of \(K\), then the natural map \(\Ext^1_K(A,B)
\to \Ext^1_L(A,B)\) need not be injective.  In fact, if \(L/K\) is a
Galois extension with Galois group \(G\), the kernel of this map can
be computed via the exact sequence
\[
0 \to H^1(G, \Hom_{L}(A,B)) \to \Ext^1_K(A,B) \to \Ext^1_L(A,B)^G.
\]
Over a number field \(K\), the group \(\Ext^1_K(A,\GG_m)\)
can be infinite whereas \(\Ext^1_K(A,B)\) is always finite.

\section{Congruences between modular forms} \label{S:ExtCong}

We will maintain the following notations throughout this section.  Let
$K$ be a number field with ring of integers $\cO$ and discriminant
$D$.  Let $A$ and $B$ be abelian varieties over $K$, and let $R>1$ be
the least integer such that both $A$ and $B$ extend to abelian schemes
(denoted $\cA$ and $\cB$) over $\cO[\frac{1}{R}]$.  Let
\[
S := R \prod_{\substack{p\ \textnormal{prime} \\
\exists \fp \mid p,\ e(\fp) \geq p-1}} p,
\]
where $e(\fp)$ is the ramification index of $\fp$ in $\cO$.
Now take $\Ext^1_{\cO[\frac{1}{S}]}(\cA, \cB)$ for the $\Ext$-group in
the category of commutative group schemes over $\cO[\frac{1}{S}]$.
The necessity of changing from $R$ to $S$ will be explained below.

\begin{proposition} \label{P:Ext1Siso}
  For abelian varieties $A$ and $B$ over a number field $K$, the
  natural map $\Ext^1_{\cO[\frac{1}{S}]}(\cA,\cB) \to \Ext^1_K(A,B)$
  is an isomorphism.
\end{proposition}

Now suppose that \(A\) and \(B\) are elliptic curves. Let $\fp$ be an
ideal of $\cO$ coprime to $S$.  Both $A$ and $B$ have good reduction
at $\fp$; let $\tA$ and $\tB$ denote the corresponding elliptic curves
over $\FF_{\fp}$. 

\begin{theorem} \label{T:ExtCong}
  Let $A$ and $B$ be elliptic curves over a number field $K$, and let
  \(e\) be the exponent of \(\Ext_K(A,B)\).  For a prime ideal $\fp$
  of $\cO$ coprime to $S$,
\[ 
\#\tA(\FF_\fp) \equiv \#\tB(\FF_\fp) \pmod{e}.
\]
\end{theorem}

This theorem has the following corollary for congruences between
Fourier coefficients of modular forms.  We provide two proofs of this
corollary.  The first is a direct application of
Theorem~\ref{T:ExtCong}.  The second proof, found by Ribet after
reading a previous version of this paper, is direct and relies on the
theta operator.

\begin{corollary} \label{C:Extmodf}  
  Suppose $A$ and $B$ are elliptic curves over $\QQ$ of conductors $M$
  and $N$ respectively.  Let $f = \sum a_n q^n$ and $g = \sum b_n q^n$
  be the associated normalized newforms of \(A\) and \(B\), and let
  $e$ be the exponent of $\Ext^1_{\QQ}(A,B)$.  Then, for any integer
  $n$ with $\gcd(n,2MN) = 1$,
\[
  a_n \equiv b_n \pmod{e}.
\]
\end{corollary} 


\begin{remark}
  It is possible to strengthen the corollary.  Since \(f\) and \(g\)
  depend only on the isogeny class of \(A\) and \(B\) and yet the
  group \(\Ext^1_{\QQ}(A,B)\) is not invariant under isogeny (in
  general), one can replace \(e\) in the corollary above by (i) the
  exponent \(e'\) of the group \(\Ext^1_{\QQ}(A',B')\) where \(A'\)
  and \(B'\) are isogenous to \(A\) and \(B\) over \(\QQ\); or, (ii)
  the least common multiple of the set of all \(e'\).
\end{remark}


\begin{proof}[Proof of Proposition~\ref{P:Ext1Siso}]
Let $n$ be any positive integer.  We first show that the natural map
\[
  b: \Hom_{\cO[\frac{1}{S}]} (\cA_n,\cB_n) \iso \Hom_K(A_n,B_n)
\]
is an isomorphism.  This follows from standard patching arguments
\cite[pp.\ 43--45]{Mazur}  once we show that, for each prime \(\fq\)
coprime to $S$, the natural map
\[
  \Hom_{\cO_{\fq}}(\cA_n, \cB_n) \iso \Hom_{K_{\fq}} (A_n,B_n)
\]
is an isomorphism.  Let $q$ be the residue characteristic of $\fq$.
The second map is clearly an isomorphism if $q \nmid n$: the \'etale
group schemes $\cA_n$ and $\cB_n$ over $\Spec \cO_{\fq}$ are
determined by the Galois modules $A_n(\Kbar_{\fq})$ and $B_n(\Kbar_{\fq})$
\cite[pp.\ 43--45]{Mazur}.  In the case that $q \mid n$, write  \(n=dr\)
with \(d\) a power of \(q\) and \(r\) coprime to \(q\). We apply
\cite[Cor.\ of Thm.~4.5.1]{tate97} to the commutative finite flat
group schemes $\cA_d$, $\cB_d$ over $\Spec \cO_{\fq}$, which implies that
the natural map
\[
  \Hom_{\cO_{\fq}} (\cA_d,\cB_d) \to \Hom_{\Gal(\Kbar_{\fq}/K_{\fq})}
  (A_d(\Kbar_{\fq}),B_d(\Kbar_{\fq}))
\]
is an isomorphism. As the result \cite[Cor.\ of Thm.~4.5.1]{tate97}
assumes that the ramification index \(e(\fq)\) of \(\cO_\fq\) is less
than $q-1$, we are forced to switch from $R$ to~$S$.

We obtain a commutative diagram from
\eqref{E:ExtABn} with exact rows
\begin{equation} \label{E:ExtStoK}
\begin{array}{c}
\xymatrix@C-15pt@R-15pt{0 \ar[r] & \Hom_{K}(A,B)\otimes {\ZZ/{n\ZZ}} \ar[r]
& \Hom_{K}(A_n, B_n) \ar[r] & 
  \Ext^1_{K}(A,B)_n \ar[r] & 0 \\
0 \ar[r] & \Hom_{\cO[\frac{1}{S}]}(\cA,\cB)\otimes
{\ZZ/{n\ZZ}} \ar[r] \ar[u]_{\wr}^{a} 
  & \Hom_{\cO[\frac{1}{S}]}(\cA_n, \cB_n) \ar[r]
\ar[u]_{\wr}^{b} & 
  \Ext^1_{\cO[\frac{1}{S}]}(\cA,\cB)_n \ar[r] \ar[u]_{\wr}^{c} & 0.}
\end{array}
\end{equation}
The second row is obtained by the analogue of \eqref{E:ExtABn} using
$0 \to \cA_n \to \cA \to \cA \to 0$ \cite[\S II.5]{milne:ADT}.  The
vertical maps are the natural restriction maps.  Now, $a$ is an
isomorphism by the N\'eron mapping properties of $\cA$ and $\cB$.  By
the preceding paragraph, $b$ is also an
isomorphism.  Therefore, $c$ is an isomorphism.
\end{proof}

\begin{proof}[Proof of Theorem~\ref{T:ExtCong}]
  Let $\fp$ be a prime of $\cO$, coprime to $S$, with residue
  characteristic $p$.  We obtain the following diagram that can be
  appended to \eqref{E:ExtStoK},
\begin{equation} \label{E:ExtRed}
\begin{array}{c}
\xymatrix@C-15pt@R-15pt
{0 \ar[r] & \Hom_{\cO[\frac{1}{S}]}(\cA,\cB)\otimes
{\ZZ/{n\ZZ}} \ar[r] \ar[d]_{h} 
  & \Hom_{\cO[\frac{1}{S}]}(\cA_n, \cB_n) \ar[r]
\ar[d]_{i} & 
  \Ext^1_{\cO[\frac{1}{S}]}(\cA,\cB)_n \ar[r]
\ar[d]_{j} & 0 \\
  0 \ar[r] & \Hom_{\FF_{\fp}}(\tA, \tB) \otimes {\ZZ/{n\ZZ}} \ar[r]
  & \Hom_{\FF_{\fp}}(\tA_n, \tB_n) \ar[r] & \Ext^1_{\FF_{\fp}}(\tA,\tB)_n
  \ar[r] & 0.} 
\end{array}
\end{equation}
Put $\alpha = h\circ a^{-1}$, $\beta = i\circ b^{-1}$, and $\gamma
=j\circ c^{-1}$.  It is important to note that
\[
\beta: \Hom_{K}(A_n,B_n) \to \Hom_{\FF_\fp}(\tA_n,\tB_n)
\]
is injective: if $\gcd(n,p) = 1$, this is clear; and if $n = p^r$,
this follows from the faithful nature of the functor $G
\rightsquigarrow G \times_{\cO_\fp} \FF_\fp$ which maps commutative
finite flat group schemes over $\cO_\fp$ of $p$-power order to their
special fiber (see \cite[Thm.~4.5]{conrad97} and \cite[Thm.~1, p.~171;
Thm.~2, p.~217]{fontaine}). These references give the required result
as the ramification index \(e(\fp)\) is less than \(p-1\).  The
injectivity of $\beta$ provides an exact sequence
\begin{equation} \label{E:CoKerSeq}
  0 \to \ker(\gamma) \to \coker(\alpha) \to \coker(\beta) \to
  \coker(\gamma) \to 0.
\end{equation}

We can now complete the proof of the theorem.  We consider the
non-trivial case that $\#\tA(\FF_\fp) \neq \#\tB(\FF_\fp)$.  This implies
that \(\Hom_{\FF_\fp}(\tA,\tB) = 0\), and so combining
\eqref{E:ExtStoK}--\eqref{E:CoKerSeq} and Theorem \ref{T:ExtNF}, we
have the injectivity of $\gamma : \Ext^1_{K}(A,B) \hookrightarrow
\Ext^1_{\FF_\fp}(\tA,\tB)$.  By Theorem~\ref{Milne}(a),
\[
  \# \Ext^1_{\FF_p}(\tA,\tB) = \bigl( \#A(\FF_{\fp}) - \#B(\FF_{\fp}) \bigr)^2.
\]

Moreover, since $A$ and $B$ are elliptic curves, Theorem
\ref{Milne}(b) and the Cassels-Tate pairing combine to give
\[
  \Sha(\tB/\FF_\fp(\tA)) \cong T \times T,
\]
for some abelian group $T$.  Thus the exponent $e$ of
$\Ext^1_{\FF_p}(\tA,\tB)$ divides $\# T = | \#A(\FF_{\fp}) -
\#B(\FF_{\fp})|$.  This, with the injectivity of $\gamma$,
gives
\[
  e \mid \bigl( \#A(\FF_{\fp}) - \#B(\FF_{\fp}) \bigr),
  \quad \forall \fp \nmid S.
\]
\end{proof}

\begin{proof}[Proof of Corollary~\ref{C:Extmodf}]
  The corollary follows almost immediately from
  Theorem~\ref{T:ExtCong}.  We need only note that since $K=\QQ$, the
  primes  which divide $S$ are exactly the primes which divide
  \(2MN\).  The congruence for all $n$ with $\gcd(n,S)=1$ 
  follows since $f$ and $g$ are Hecke eigenforms.
\end{proof}

\begin{remark}
  Suppose that $K= \QQ$ and that \(R\) is odd. (i) The map
  $\Ext^1_{\ZZ[\frac{1}{R}]}(\cA,\cB) \to \Ext^1_{\QQ}(A,B)$ may not
  be an isomorphism.  However, it is injective. As \(a\) is an
  isomorphism even over $\ZZ[\frac{1}{R}]$ by the N\'eron mapping
  property, it suffices to show that the map
  $\Hom_{\ZZ[\frac{1}{R}]}(\cA_n, \cB_n) \to \Hom_{\QQ}(A_n,B_n)$ is
  injective.  Given that $b$ is injective over $\ZZ[\frac{1}{S}]$, it
  suffices to check that the map $\Hom_{\ZZ_2}(\cA_n,\cB_n) \to
  \Hom_{\QQ_2}(A_n,B_n)$, with $n = 2^r$, is injective.  This is clear
  \cite[p.~152]{tate97}.  Thus, the $\Ext$-group over
  $\ZZ[\frac{1}{R}]$ might be smaller than
  $\Ext^1_{\cO[\frac{1}{S}]}(\cA,\cB)$, and so a potential analogue of
  Theorem~\ref{T:ExtCong} will be weaker.
  
  (ii) In addition, the definition of $\beta$ over $\ZZ[\frac{1}{R}]$
  is problematic, due to the non-exactness of the N\'eron model
  functor for $p=2$ (see~\cite[Ex.~4, p.~190]{blr} for a
  counterexample to \cite[Thm.~4, p.~187]{blr}).
\end{remark}

\begin{proof}[Direct Proof of Corollary~\ref{C:Extmodf}] (This
  proof is due to Ribet.)
Since $f$ and $g$ are Hecke eigenforms, it suffices to show
that $a_p \equiv b_p \pmod{\lambda}$ for all primes $p \nmid 2MN$ and
prime powers $\lambda = \ell^m$ dividing $e$.  Fix now such $p$ and
$\lambda = \ell^m$.

Suppose first that $p \neq \ell$.  Let $G_{\QQ} := \Gal(\Qbar/\QQ)$.
Now letting $\sigma \in G_{\QQ}$ be a Frobenius
element at $p$, the Eichler-Shimura relations
\cite[Thm.~4.2.2]{hida:GMF} dictate that
\begin{align*}
  \sigma^2 - a_p \sigma + p &= 0 \quad \textnormal{on $A_\lambda$,} \\
  \sigma^2 - b_p \sigma + p &= 0 \quad \textnormal{on $B_\lambda$.}
\end{align*}
Now from \eqref{E:ExtABn} with $n = \lambda$, it follows that there is
a Galois equivariant morphism $\alpha : A_\lambda \to B_\lambda$ with
exponent $\lambda$.  Because $\alpha$ is $G_{\QQ}$-equivariant, it
follows that the operators $a_p \sigma$ and $b_p \sigma$ agree on the
image of $\alpha$ in $B_\lambda$, which has exponent $\lambda$.  Hence
$a_p \equiv b_p \pmod{\lambda}$.

Now suppose that $p = \ell$.  In particular, $\ell \nmid 2MN$.
Now by the preceding paragraph, we see that $a_n \equiv b_n
\pmod{\lambda}$ whenever $n$ is coprime to $2MN\ell$.  Let $\af$ and $\ag$
be modular forms with Fourier expansions
\[
  \af = \sum_{\substack{n=1 \\ \gcd(n,2MN)=1}}^{\infty} a_n q^n,
\qquad
  \ag = \sum_{\substack{n=1 \\ \gcd(n,2MN)=1}}^{\infty} b_n q^n.
\]
(These modular forms can be obtained by twisting twice by appropriate
quadratic characters, and as such have level dividing
$4M^2N^2$.)  Thus, if we let $\Theta = q \frac{d}{dq}$ be the usual
$\Theta$ operator on modular forms, we see that $\Theta^m$ annihilates
$\af-\ag$ modulo~$\lambda$.  In particular, as power series in $q$,
\[
  \Theta(\af - \ag) \equiv 0 \pmod{\ell}.
\]
Since $\ell$ is odd, the $\Theta$ operator is injective modulo $\ell$
on modular forms of weight~$2$ \cite[\S II]{katz77}.  Therefore $\af
\equiv \ag \pmod{\ell}$.  Moreover, if $m > 1$, we see that $\Theta$
annihilates $(\af - \ag)/\ell$ modulo $\ell$, and thus $\af \equiv \ag
\pmod{\ell^2}$.  By induction, $\af \equiv \ag \pmod{\lambda}$.
\end{proof}

\begin{remark} (Ribet) When $e$ is odd, the
arguments in both proofs of Corollary~\ref{C:Extmodf} imply that $a_n
\equiv b_n \pmod{e}$ for $\gcd(n,MN) = 1$ instead of $\gcd(n,2MN) =
1$. In
some sense, having the $\Theta$ operator injective modulo odd primes
is similar to having the absolute ramification index of $\ZZ_p$
strictly less than $p-1$ in the first proof of Corollary~\ref{C:Extmodf}.
\end{remark}

\section{Modular parametrizations and congruence moduli} \label{S:ParMod}

We now investigate how the theory of congruence moduli of Hida
\cite{hida81}, Ribet \cite{ribet84} and D.~Zagier \cite{zag85} fits in
with the \(\Ext\)-group.

Let $J_0(N)$ be the Jacobian of the modular curve $X_0(N)$ over $\QQ$,
and let $S_2(\Gamma_0(N))$ be the space of weight two cusp forms on
$\Gamma_0(N)$.  For an elliptic curve $A$ over $\QQ$ of conductor $N$,
let $f = \sum a_n q^n \in S_2(\Gamma_0(N))$ be its associated
normalized newform.  The \emph{congruence modulus of $A$} \cite[\S
5]{zag85} is
\[
m_A := \max\left\{m \in \ZZ \biggm|
\begin{tabular}{c}
$\textstyle \exists g = \sum b_n q^n\in (f)^{\perp} \cap
\power{\ZZ}{q}$ so that \\[2pt]
$b_n \equiv a_n \pmod{m}$ for all $n$
\end{tabular}
\right\},
\]
where $(f)^{\perp} \subseteq S_2(\Gamma_0(N))$ is the orthogonal
complement with respect to the Petersson inner product.  The
\emph{restricted congruence modulus of $A$} is
\[
r_A := \max\left\{r \in \ZZ \biggm| 
\begin{tabular}{c}
$\textstyle \exists g = \sum b_n q^n\in (f)^{\perp} \cap
    \power{\ZZ}{q}$ so that \\[2pt]
$b_n \equiv a_n \pmod{r}$ for all $n$ with $\gcd(n,2N) = 1$
\end{tabular}
\right\}.
\]

Assume that $A$ is an optimal quotient of $J_0(N)$ (i.e.\ a strong
Weil curve) and that $\phi_A : X_0(N) \to A$ is its modular
parametrization. This induces an exact sequence of abelian varieties
over~$\QQ$,
\begin{equation}\label{eta}
  \eta : \quad  0 \to C \to J_0(N) \stackrel{\phi_*}{\to} A \to 0.
\end{equation}
A well-known result of Ribet and Zagier \cite{zag85} is that:
\emph{the degree $d_A$ of $\phi_A$ divides \(m_A\)}, i.\ e.\ 
\[
  d_A \mid m_A.
\]
Of course $m_A$ divides $r_A$.  Our next result, that the exponent
$e_{A}$ of $\Ext^1_{\QQ}(A,C)$ sits between $d_A$ and $r_A$,
ultimately relies on the Eichler-Shimura relations.

\begin{theorem} \label{T:ParExtRestMod}
With notations as above, one has 
\[
  d_A \mid e_A \mid r_A.
\]
\end{theorem}

\begin{remark}Since $d_A$ is unbounded, 
Theorem~\ref{T:ParExtRestMod} shows that the
order of $\Ext^1_{\QQ}(A,C)$ is unbounded.
\end{remark}

In general, $e_A$ does not divide the congruence modulus $m_A$.  The
analogous result for $m_A$ requires a refinement of
$\Ext^1_{\QQ}(A,C)$, which involves the Hecke algebra $\TT = \ZZ[T_1,
T_2, \ldots ]$ associated with $J_0(N)$.  Namely, let
$\Ext^1_{\QQ,\TT}(A,C)$ be the Yoneda $\Ext$-group in the category of
commutative $\TT$-group schemes over $\QQ$, i.e.\ groups $G$ together
with a homomorphism $\TT \to\End_{\QQ}(G)$. Write $e_{A,\TT}$ for the
exponent of $\Ext^1_{\QQ,\TT}(A,C)$. We show below that the natural
forgetful map $\Ext^1_{\QQ,\TT}(A,C)\to \Ext^1_{\QQ}(A,C)$ is
injective, and thus $e_{A,\TT} \mid e_A$.

\begin{theorem} \label{T:ParExtTMod} With notations as above, one has 
\[
  d_A \mid e_{A,\TT} \mid m_A.
\]
\end{theorem}

The sequence (\ref{eta}) is also an exact sequence of $\TT$-group schemes
over~$\QQ$.  For any integer~$n$, we have the following commutative
diagram with exact rows,
\[
\xymatrix@C-14pt@R-14pt{0 \ar[r] & \Hom_{\QQ,\TT}(A,C)\otimes
  {\ZZ/{n\ZZ}} \ar[r] \ar[d] & \Hom_{\QQ,\TT}(A_n, C_n) \ar[r]
  \ar[d]_{\kappa} &
  \Ext^1_{\QQ,\TT}(A,C)_n \ar[r] \ar[d]_{\iota} & 0 \\
  0 \ar[r] & \Hom_{\QQ}(A, C) \otimes {\ZZ/{n\ZZ}} \ar[r] &
  \Hom_{\QQ}(A_n, C_n) \ar[r] & \Ext^1_{\QQ}(A,C)_n \ar[r] & 0.}
\]
The bottom row is precisely \eqref{E:ExtABn}, and since the Kummer
sequence $0 \to A_n \to A \to A \to 0$ is a sequence of $\TT$-group
schemes, the top row is obtained in the same manner
as~\eqref{E:ExtABn} but with the requirement of $\TT$-equivariance.
Note also that the top row is an exact sequence of $\TT$-modules. Now
since $\Hom_{\QQ}(A,C) = 0 = \Hom_{\QQ,\TT}(A,C)$ and the map $\kappa$
is injective, it follows that $\iota$ is also injective.
Alternately, any element of $\Ext^1_{\QQ,\TT}(A,C)$ which splits in
$\Ext^1_{\QQ}(A,C)$ is split by a $\TT$-equivariant morphism.  By
Theorem \ref{T:ExtNF}, the injectivity of $\iota$ shows that
$\Ext^1_{\QQ,\TT}(A,C)$ is finite.

\begin{proposition} \label{P:ParExt}
  The order of $\eta$ in $\Ext^1_{\QQ,\TT}(A,C)$ is equal to the
  degree of $\phi_A$.
\end{proposition}

\begin{proof}
  Let $d_A = \deg \phi_A$.  Because $A$ is an optimal quotient of
  $X_0(N)$, the dual $\phi^*: A \to J_0(N)$ of $\phi_*$ is injective,
  and the composition $\phi_* \circ \phi^* \in \Hom_{\QQ}(A,A)$ is
  multiplication by $d_A$.  From the exact sequence $0 \to A \to
  J_0(N) \to C^{*} \to 0$ and the fact that $\Hom_{\QQ}(A,C^{*}) = 0$,
  it follows that the map
\[
\Hom_{\QQ}(A,A) \iso \Hom_{\QQ}(A,J_0(N))
\]
induced by $\phi^*$ is an isomorphism.  Furthermore, the map
\[
\Hom_{\QQ}(A,A) \to \Hom_{\QQ}(A,J_0(N)) \to \Hom_{\QQ}(A,A),
\]
induced by $\phi_* \circ \phi^*$, is multiplication by $d_A$, and so
the image of $\Hom_{\QQ}(A,J_0(N))$ is precisely
$d_A(\Hom_{\QQ}(A,A))$.

The long exact sequence for $\Ext^i_{\QQ}(A,-)$ applied to $\eta$ begins
\[
0 \to \Hom_{\QQ}(A,J_0(N)) \stackrel{\phi_*}{\to} \Hom_{\QQ}(A,A) \to
\Ext^1_{\QQ}(A,C),
\]
and the image of an endomorphism $\alpha \in \Hom_{\QQ}(A,A)$ in
$\Ext^1_{\QQ}(A,C)$ is the pull-back $\alpha^{\ast}\eta$.  Therefore,
because $\Hom_{\QQ}(A,A) = \ZZ$, it follows that the image of
$\Hom_{\QQ}(A,A)$ in $\Ext^1_{\QQ}(A,C)$ is the subgroup generated by
$\eta$.  Thus the order of $\eta$ in $\Ext^1_{\QQ}(A,C)$ is $d_A$, and
since $\eta$ represents a class in $\Ext^1_{\QQ,\TT}(A,C) \subseteq
\Ext^1_{\QQ}(A,C)$, we are done.
\end{proof}

Now there are two ways to define a $\TT$-module structure on
$\Ext^1_{\QQ}(A,C)$, namely by pushing out along $C$ or pulling back
along $A$.  These two $\TT$-module structures need not coincide.
However, because $\Ext^1_{\QQ,\TT}(A,C)$ consists of classes of
extensions that are $\TT$-equivariant, it follows that the two
$\TT$-module definitions restricted to $\Ext^1_{\QQ,\TT}(A,C)$
\emph{are} the same.  In fact, $\Ext^1_{\QQ,\TT}(A,C)$ is the largest
subgroup of $\Ext^1_{\QQ}(A,C)$ on which the two $\TT$-module
structures agree.

\begin{proof}[Proof of Theorem~\ref{T:ParExtTMod}]
By Proposition~\ref{P:ParExt} it suffices now to prove that $e_{A,\TT} \mid
m_A$.  Let $I_A$ and $I_C$ be the kernels of the maps $\TT \to
\End_{\QQ}(A)$ and $\TT \to \End_{\QQ}(C)$.  By the general
considerations of \cite[\S 2]{diam89}, we see that $m_A$ is the
exponent of $\TT_{A,C} := \TT/(I_A + I_C)$ as an abelian group.  Now
clearly both $I_A$ and $I_C$ annihilate $\Ext^1_{\QQ,\TT}(A,C)$, and
so $\Ext^1_{\QQ,\TT}(A,C)$ is a $\TT_{A,C}$-module.  Thus we must have
$e_{A,\TT} \mid m_A$.
\end{proof}

\begin{proof}[Proof of Theorem~\ref{T:ParExtRestMod}]
Again by Proposition~\ref{P:ParExt} it suffices to prove that $e_A
\mid r_A$.  Let $\TT' \subseteq \TT$ be the subalgebra generated
by all $T_n$ with $\gcd(n,2N) = 1$, and let $I_A'$ and $I_C'$
be the kernels of the maps $\TT' \to \End_{\QQ}(A)$ and $\TT'
\to \End_{\QQ}(C)$.  From \cite[\S 2]{diam89}, it follows that $r_A$
is the exponent of $\TT'_{A,C} := \TT'/(I_A' +
I_C')$.  If we can show that $\Ext^1_{\QQ}(A,C)$ is a
$\TT'_{A,C}$-module, then we are done.  Namely, we need to show
that the two operations of $\TT'$ on $\Ext^1_{\QQ}(A,C)$
coincide.

Consider $T_{p}$ with $p \nmid 2N$.  Since $\Hom_{\QQ}(A,C) = 0$,
\eqref{E:ExtABn} implies that $\Ext^1_{\QQ}(A,C)_n \cong
\Hom_{\QQ}(A_n,C_n)$.    From the
Eichler-Shimura relations \cite[Thm.~4.2.1]{hida:GMF}, we have that
\begin{equation} \label{E:EichShim}
  T_{p} = F + V \in \End_{\FF_{p}}(\tJ),
\end{equation}
where $\tJ$ is the reduction of $J_0(N)$ modulo $p$, and $F$ and $V$
are the $p$-th power Frobenius and the Verschiebung on $\tJ$.  As in
\eqref{E:ExtRed}, the natural map
\begin{equation} \label{E:RedInj}
  \Hom_{\QQ}(A_n,C_n) \hookrightarrow \Hom_{\FF_{p}}(\tA_n,\tC_n)
\end{equation}
is injective.  Since both $A$ and $C$ are subabelian varieties of
$J_0(N)$, the action of $T_p$ on $\tA_n$ and $\tC_n$ is determined by
the restriction of \eqref{E:EichShim}. In particular, for any $\alpha
\in \Hom_{\QQ}(A_n,C_n)$, its image $\tilde{\alpha} \in
\Hom_{\FF_{p}}(\tA_n,\tC_n)$ satisfies
\[
\tilde{\alpha} \circ T_{p}|_{\tA_n} = T_{p}|_{\tC_n} \circ
\tilde{\alpha}.
\]
The injectivity of \eqref{E:RedInj} shows that in fact $\alpha \circ
T_{p}|_{A_n} = T_{p}|_{C_n} \circ \alpha$, and we are done.
\end{proof}

\begin{corollary}
  If the conductor $N$ of $A$ is square-free, then $d_A = e_{A,\TT}=
  m_A$.
\end{corollary}

\begin{proof}
  For prime $N$, Ribet and Zagier \cite[Thm.~3]{zag85} have shown that
  $d_A = m_A$.  For square-free $N$, that $d_A=m_A$ has been proved
  recently by Ribet \cite{ribetnew} as a consequence of his proof of a
  conjecture of A.~Agashe and W.~Stein.
\end{proof}

\begin{remark}
  If $N$ is prime and if the form $g = \sum b_n q^n\in (f)^{\perp}
  \cap \power{\ZZ}{q}$, which gives $b_n \equiv a_n \pmod{r_A}$ for all
  $n$ with $\gcd(n,2N) = 1$, also satisfies $b_2 \equiv a_2
  \pmod{r_A}$, then in fact $d_A = e_{A,\TT} = m_A = e_A = r_A$ by the
  Sturm bound \cite{sturm87}.
\end{remark}
 
\begin{example}
Theorems~\ref{T:ParExtRestMod} and~\ref{T:ParExtTMod} are best
possible in the following sense.  The mod~$3$ representation of the
elliptic curve $90C1$ in \cite{crem97} of conductor $90$,
\[
  A : y^2 + xy + y = x^3 - x^2 +13x - 61,
\]
has image $\left( \begin{smallmatrix} * & * \\
    0 & 1 \end{smallmatrix} \right)$ \cite[Thm.~3.2]{revvil01}.
Likewise, the mod~$3$ representation of the elliptic curve $90A1$ in
\cite{crem97}, also of conductor $90$,
\[
  B : y^2 + xy = x^3 - x^2 + 6x,
\]
has image $\left( \begin{smallmatrix} 1 & * \\ 0 & *
\end{smallmatrix} \right)$.
Both $A$ and $B$ are strong Weil curves. By their mod $3$
representations, there is a non-trivial element $\alpha \in
\Hom_{\QQ}(A_3,B_3)$, which induces a homomorphism $\alpha : A_3 \to
C_3$, where $C$ is the kernel of $J_0(90) \to A$.  Thus
$\Ext^1_{\QQ}(A,C)_3 \cong \Hom_{\QQ}(A_3,C_3)$ has a non-trivial
element of order $3$.  However, direct computation provides $d_A = 16$
and $m_A = 16$. In particular, $m_A \neq r_A$ and $e_A\nmid m_A$.
\end{example}

\begin{remark}
If $A$ is an elliptic curve over $\QQ$ and 
$\phi:J \to A$ realizes $A$ as an optimal quotient of the
  Jacobian $J$ of a Shimura curve, then we have an exact sequence 
\[ 
\eta : \quad 0 \to C \to J \stackrel{\phi}{\to} A \to 0.
\]
Theorem~\ref{T:ParExtRestMod} suggests that the exponent of
$\Ext^1_{\QQ}(A,C)$ provides a ``congruence modulus''.  In fact, the
Hida constant $c_A$ defined by E. Ullmo \cite[p.~326]{ullmo01} is the
order of $\eta \in 
\Ext^1_{\QQ}(A,C)$.
\end{remark}

\section{Special values of \(L\)-functions}\label{BlochKato}

We provide a restatement in terms of \(\Ext\)-groups of the deep
results of Diamond, Flach, and Guo \cite{dfg01, dfg} on the Bloch-Kato
conjecture.

\subsection*{Symmetric square of an elliptic curve}

Let $E$ be an elliptic curve over $\QQ$ of conductor $N$ with
$\End_{\Qbar}(E) = \ZZ$. Consider the symmetric square $L$-function
$L(\Sym^2 E, s)$ of $E$ \cite{flach92}.  Note that
\[
  L(\Sym^2E,s) = L(\Sym^2h_1(E),s) = L(\Sym^2h^1(E), s+2).
\]
The Bloch-Kato conjecture
\cite[Eq.~(2)]{flach92} on the special value at \(s=2\) of $L(\Sym^2
E, s)$ states that
\begin{equation} \label{E:flachL2}
\frac{L(\Sym^2 E,2)}{\Omega(2)} = \frac{
  \#\Sha({\QQ},A(2))}{
  \#H^0({\QQ},A (1))\cdot \#H^0(\QQ,A(2))} 
  \prod_p c_p.
\end{equation}
For the sake of brevity we refer to Flach \cite[\S\S 0--1]{flach92}
for definitions, and we use his notation.

For a field $K$, set $G_{K} := \Gal(\Kbar/{K})$.  Let $P_E$ denote the
finite set of rational primes consisting of (a) all $\ell \mid 2N$;
and (b) all $\ell$ such that the Galois representation on $E_\ell$
restricted to $G_F$, where $F = \QQ(\sqrt{(-1)^{(\ell-1)/2}\,\ell})$,
is not absolutely irreducible.  For $\ell \notin P_E$, the $\ell$-part
of the Bloch-Kato conjecture has been proved by Diamond, Flach, and
Guo \cite[Thm.~0.2]{dfg01}, \cite[Thm.~8.9]{dfg}.  The following
theorem is a reformulation of their result in terms of
\(\Ext\)-groups.

\begin{theorem} \label{T:Sym2}
  If $E$ is an elliptic curve over \(\QQ\) with \(\End_{\Qbar}(E) =
  \ZZ\), then up to powers of $\ell$ for $\ell \in P_E$, 
\begin{equation} \label{E:ExtL2}
\frac{L(\Sym^2 E,2)}{\Omega(2)} = \frac{
  \#\Sha(\Ext^1_{\Qbar}(E,E))}{
  \#\Ext^1_{\QQ}(E,E)\cdot \#\Ext^1_{\Qbar}(E,E)(1)^{G_{\QQ}}} 
  \prod_p c_p.
\end{equation}
\end{theorem}

\begin{proof}
  Because of the results \cite[Thm.~0.2]{dfg01} and
  \cite[Thm.~8.9]{dfg}, it suffices to match the terms in
  \eqref{E:flachL2} with those of \eqref{E:ExtL2}.
  
  The total Tate module $TE$ of \(E\) is a rank two module over
  \(\widehat{\ZZ}:= \varprojlim {\ZZ}/{n\ZZ}\).  By the Weil pairing,
  \(TE\) is isomorphic to the Tate twist $TE^{\vee}(1)$ of the dual
  module $TE^{\vee}$.  From \eqref{E:TpExt}, we see that
  $T\Ext^1_{\Qbar}(E,E)$ is isomorphic to the quotient of
  \(\End_{\widehat{\ZZ}}(TE)\) by the \(\widehat{\ZZ}\)-submodule
  generated by the identity map.  Via the self-duality
  \(\End_{\widehat{\ZZ}}(TE)\cong TE\otimes TE^{\vee}\), we have
\begin{equation} \label{E:TExt}
  T \Ext^1_{\Qbar}(E,E) \cong (\Sym^2 TE)(-1).
\end{equation}
The terms in the denominator of \eqref{E:flachL2} can easily be
identified as 
\begin{align*}
\#H^0({\QQ},A(1)) &= \# H^0(G_{\QQ}, (\Sym^2 TE)(-1) \otimes \QQ/\ZZ),\\
\#H^0(\QQ,A(2)) &= \# H^0(G_{\QQ}, (\Sym^2 TE)\otimes \QQ/\ZZ).
\end{align*}
By \eqref{E:TExt}, the former is $\# \Ext^1_{\Qbar}(E,E)^{G_\QQ}$, and
the latter is $\# \Ext^1_{\Qbar}(E,E)(1)^{G_\QQ}$.  Moreover, the
numerator of \eqref{E:flachL2} is
\[
\# \Sha((\Sym^2 TE)\otimes \QQ/\ZZ) = \# \Sha((\Sym^2 TE)(-1)
\otimes \QQ/\ZZ),
\]
where the equality follows from \cite[Thm.~1]{flach90}.  The proof of
the theorem is then complete by the following lemma.
\end{proof}

\begin{lemma} \label{L:ExtInv}
  Let $E$ be an elliptic curve over a number field $K$ with
  $\End_{\Kbar}(E) = \ZZ$.  For any odd integer $n$, $\Ext^1_K(E,E)_n
  \cong \Ext^1_{\Kbar}(E,E)_{n}^{G_K}$.
\end{lemma}

\begin{proof}
By taking the long exact sequence of $G_K$-cohomology of
  \eqref{E:ExtABn} with $K = \Kbar$, we obtain
  an exact sequence
\begin{multline*}
  0 \to \Hom_K(E,E) \otimes_{\ZZ} \ZZ/n\ZZ \to \Hom_K(E_n,E_n)
    \to \Ext^1_{\Kbar}(E,E)^{G_K}_{n} \\
  \to \Hom(G_K,\ZZ/n\ZZ) \stackrel{\phi}{\to}
    H^1(G_K,\Hom_{\Kbar}(E_n,E_n)).
\end{multline*}
If $f \in \ker \phi$, then $\phi(f)$ is represented by a coboundary
which takes values in the trace~$0$ space of $\Hom_{\Kbar}(E_n,E_n)$.
However, by definition $\phi(f)(\sigma) =
f(\sigma)\cdot\operatorname{id}$, for $\sigma \in G_K$.  Thus for all
$\sigma \in G_K$, we have $2f(\sigma) = 0$, which implies $f=0$.
Therefore, $\phi$ is injective, and the result follows from
\eqref{E:ExtABn}.
\end{proof}

\end{document}